\documentclass{article}%
\usepackage{array}
\usepackage{hyperref}
 \hypersetup{colorlinks   = true, urlcolor= blue, linkcolor=black,citecolor=blue}
\usepackage{graphicx}
\usepackage[numbers]{natbib}
\usepackage{amsmath}
\usepackage{amssymb} 
\usepackage{bm} 
\usepackage{geometry} 
\usepackage[section]{placeins} 
\usepackage{caption}
\usepackage[singlelinecheck=on,labelformat=simple]{subcaption}

\begin{document}
\title{Effect of Transit Signal Priority on Bus Service Reliability}
\author{Paul Anderson, anderson@berkeley.edu\footnote{Corresponding author}\\
\normalsize University of California, Berkeley\\
Carlos F.  Daganzo, daganzo@berkeley.edu\\
\normalsize University of California, Berkeley\\
}
\maketitle

\section*{Abstract}

Buses are known to be an unstable system, and as a result they struggle to stay on schedule and/or to maintain regular headways. In current practice, transit agencies add significant slack to the schedule and hold buses until their scheduled departure time at certain stations along the route. These practices can smooth out small disruptions but are not robust to large disruptions. Various headway-based strategies have been proposed but these also rely on holding buses and therefore reduce their average pace. Transit signal priority (TSP) is commonly used to reduce the signal delay that buses experience but the potential for signals to serve as an additional control agent to enhance bus reliability has not been systematically evaluated. We are especially interested in evaluating conditional signal priority (CSP), in which buses send priority requests only when a condition is met. A mathematical model based on Brownian motion is proposed for bus systems with a schedule where pairing is not a concern and is used to develop formulas for the expectation and variance of steady state deviation. Simulation is used to confirm these findings and to evaluate the case of bus systems operated by headway where pairing is a concern. The analytical and simulation results show that CSP can be used not just to accelerate buses but to improve their reliability. CSP is as good as TSP in some scenarios and better in others. In all scenarios it is better for other traffic because it sends fewer priority requests.

\indent \textbf{Keywords:} public transit; bus operations; signal priority; reliability

\section{Introduction}

Buses are known to be an unstable system. Even small variations in traffic and station dwell times can be enough to throw buses off schedule or disrupt their headways. Once perturbed, they struggle to recover without, and sometimes even with, active intervention. For bus operators, a consequence of instability is that they must plan for bus runs of longer duration, which increases the cost of operating a given service. For passengers, riding times become less predictable, and the expectation and variance of out-of-vehicle waiting times also increase. The effect on waiting is particularly undesirable because this is the part of the trip that passengers find most onerous.

In current practice, transit agencies fight instability by publishing schedules that include considerable slack, and enforcing adherence by holding early buses at designated bus stops (stations) called control points \citep{newell1977unstable}. In sufficient amounts the slack allows most buses to arrive slightly early so they can depart the control points on time. Unfortunately, slack makes the schedule slower -- in the sense that buses are forced to travel more slowly than they can. This reduces bus productivity and increases passenger riding times, so slack can only be provided sparingly. As a result, the strategy often fails -- the strategy also fails if a major disruption such as a snow storm severely delays all buses as then the schedule becomes useless. Given the schedule's inability to work under all conditions, methods that dispose with it have been proposed.  Examples include  \citep{daganzo2009headway}, which proposed a formula for determining holding times based only on the forward headway (the headway between the current bus and the bus in front), \citep{ bartholdi2012self} which considered only the backward headway (the headway between the current bus and the bus behind), and \citep{daganzo2011reducing, xuan2011dynamic, daganzo2017unpublished}, which considered both the forward and backward headways. Some versions of these strategies have been shown to be robust to large disruptions and succeeded in field trials \citep{argote2015dynamic}, but have not yet progressed to large-scale implementation. Unfortunately, there is no way of speeding up late buses with bus actions alone; so even the most successful strategies are based on holding buses, and this reduces their productivity.

In view of all this, it makes sense to use transit signal priority (TSP) as an additional control agent to speed up late buses and in this way enhance reliability. Of course, TSP has been around for a long time, and considerable research has been devoted to the subject, but most works only focus on TSP's impact on traffic and bus speeds. Very few examine its impact on bus service reliability, and none appear to have studied the issue systematically. Therefore, this paper will attempt to fill the gap.

To minimize TSP's effect on traffic, \citep{janos2002bus} proposed to trigger the TSP priority requests only when buses were late. Since this conditional form of TSP (CSP for short) can be implemented with today's technology, and has no real downside, it will be the form of TSP considered in this paper.  Reference \citep{janos2002bus}, however, did not focus on reliability and was only a proof-of-concept: it examined how the schedule should change to accommodate the faster speeds achieved with CSP, and did so only for a particular corridor in San Juan, Puerto Rico. More recently, two other case studies \citep{ma2010development,chow2017multi} have considered reliability. Both assume that CSP is used not just to request priority but also to artificially delay early buses, e.g. by extending red phases; and also assume that this is the only way to retard buses, with no holding at stations. The approach in these two works seems ambitious; especially since both methods require considerable signal intelligence: \citep{ma2010development} uses a predictive multi-signal strategy that requires future knowledge; and \citep{chow2017multi} a non-predictive strategy that attempts to emulate the forward-headway holding control policy proposed in \citep{daganzo2011reducing}. The methods also seem somewhat impractical because signals may harm car traffic when retarding buses\footnote{Reference  \citep{anderson2017conditional} showed that delay requests are far more disruptive than priority requests and that two strategies that used them were worse for other traffic than TSP, even when their activation rate was much lower (e.g. 28\%).}. Since this collateral damage could be avoided with the common practice of holding buses at stations, the form of CSP used in this paper will use only priority requests as in \citep{janos2002bus}, combined with station holding.  The key question will be: how much reliability is gained by adding CSP to a bus route that is currently being controlled by holding at stations? 

A family of CSP policies will be presented that work with any type of holding control. Performance formulas will be developed with analysis and simulations. These formulas will be used to optimize the CSP policies and evaluate their marginal benefit under various forms of holding.

\section{The CSP recipe}

To maximize widespread use, we propose a scalable form of CSP that can be implemented on bus and street networks of any size, using any holding recipe. The technological requirements are modest and could be met by existing TSP hardware and automatic vehicle location (AVL) platforms. Reasonably, we assume that buses can communicate their locations through this AVL platform to a central controller in real time with negligible latency, and transmit one-way binary information to the signals. Most of the CSP intelligence shall reside in the buses, so signals need only respond to priority requests. To avoid scaling issues the proposed CSP policies will be logically simple, use local information, be computationally distributed and transmit minimal data. 

The best way to coordinate the actions of drivers and signals is to have a shared piece of information that is used by both agents to make control decisions. This shall be the output, $D_{n,x}$,  of a function that can be evaluated by any bus $n$ a little before reaching any location $x$ along its route. The piece of information shall represent an amount of delay that, imparted to bus $n$ at $x$, would greatly benefit the system. This desired delay is only a target and therefore is allowed to be positive or negative. Of course, $D_{n,x}$ shall only be evaluated when needed for control purposes; i.e., when bus $n$ is approaching either a signal, $i$, (so $x = i$)\footnote{An accurate estimate of the bus' arrival time at the signal is necessary to avoid unnecessary disruption. Priority requests should be sent from an upstream point at which the remaining travel time to the signal is relatively deterministic. A suitable location might be the nearest upstream station or the nearest upstream signal, whichever is closer.} or a station, $s$, where the bus can be held (so $x = s$)\footnote{Stations where control can be applied will also be called ``control points''.}. 

The desired delay is calculated in real-time as a function of dynamically available information. For bus $n$ this shall be the history of previous arrivals at $x$,  \{..., $a_{n-1,x}$\} and the series of projected arrivals \{$\alpha_{n,x}, ...$\}. We also allow the function to be space-dependent so that:

\begin{equation}\label{one}
D_{n,x}=F_x(..., a_{n-1,x}, \alpha_{n,x}, \alpha_{n+1,x}, ...).
\end{equation}

\noindent This function will be used to determine both, bus holding times at control points and CSP requests at signals. 

Since the holding times cannot be negative, they are calculated by truncating the desired delay. Thus, the holding time, $H_{n,s}$ at a control point $x = s$, is:

\begin{equation}\label{two}
H_{n,s}=[D_{n,s}]^+.
\end{equation}

\noindent The recipe defined by (1) and (2) is very general, and encompasses most if not all the holding strategies proposed to date. For example, the special case with no control is obtained by choosing $F_s$ so  $D_{n,s} \equiv 0$. And if a schedule $\{t_{n,s}\}$ is the basis for holding, then by choosing:

\begin{equation}
D_{n,s}=F_s(\alpha_{n,s})=t_{n,s}-\alpha_{n,s} \equiv -\varepsilon_{n,s}.
\end{equation}

\noindent The symbol $\varepsilon$ represents the bus lateness relative to the schedule. It will be used repeatedly in the paper.

Now consider signals and the CSP requests. In the spirit of minimal signal intelligence and communication requirements, the target $D_{n,i}$ will only be used to decide whether to trigger a request for priority.  The specific binary rule is as follows:

\begin{equation}\label{4a}
\{\text{request priority if and only if:} \hskip1mm D_{n,i} < -\delta\},
\tag{4a} 
\end{equation}

\noindent or, if a schedule is used:

\begin{equation}\label{4b}
\{\text{request priority if and only if:} \hskip1mm  \varepsilon_{n,i} > \delta \}.
\tag{4b} 
\end{equation}
\setcounter{equation}{4} 

\noindent In these expressions $\delta$ is a control parameter that can be freely chosen. Note from (4) that if we set $\delta = \infty$ priority is never requested. This setting can be used to model the basic case with no transit priority (NTP), which will be used as a benchmark. At the other extreme, if we set $\delta = - \infty$, priority is always requested. This second benchmark will be called ``TSP''.

\section{Analysis}

The paper will analyze CSP under three bus holding strategies: (i) no holding; (ii) holding by schedule; and (iii) holding by headways without a schedule. Case (i) is operationally the worst-performing, and for this reason the one where CSP offers the most promise; the next section shows how, indeed, CSP greatly improves performance in this case. Case (ii) is the most common but not robust to large disruptions. Accordingly, the paper will evaluate CSP under small disruptions. Case (iii) is appealing, especially for systems with short headways, as it provides a robust anti-bunching performance even under large disruptions. For this case, the paper will show how CSP improves both the buses' commercial speed and the headway variance, considering large disruptions.

The paper's results will include analytical formulas based on Brownian motion theory, and the results of simulations enriched with dimensional analysis. The CSP strategy will also be optimized, and then compared to the NTP and TSP benchmarks.  Below we do this for case (i), illustrating that stability can be achieved without holding at stations or delaying buses with signals.

\section{Case (i): No holding at stations}

We assume here that there is a schedule, and that drivers only follow it when dispatched from the origin. This situation could describe the behavior of buses when traveling between consecutive control points, if the control points are far apart and buses do not have any cruising guidance. We also assume that demand is light and/or the headways long so the bunching tendency can be neglected.  This allows us to model buses in isolation. Therefore, the subscript $n$ will be dropped in this section.  

We shall also assume for evaluation purposes that there is one signal per station and that the traffic signals encountered along the route have random offsets. This is reasonable because our buses cannot normally benefit from progression schemes due to their stops at stations. The assumption is also useful because it allows us to model a bus trip as a process with independent increments; i.e., where the times to go from one signal to the next, including stops, are mutually independent random variables. We shall finally assume that the bus route is homogeneous, with evenly spaced signals and other characteristics, so the random variables are also identically distributed. We choose the separation between signals as the unit of distance.

The following subsection, 4.1, introduces a mathematical model based on Brownian Motion and derives formulas for bus performance in terms of the control parameters. Subsection 4.2 tests the analytical findings of subsection 4.1 in simulation. Subsection 4.3 discusses the tradeoffs inherent in choosing the control parameters.

\subsection{Brownian Motion}

It should perhaps be intuitive that under the stipulated conditions the time-space trajectory of a bus can then be idealized as the realization of a Brownian motion process with time as the state and distance $x$ as the parameter; see \citep{newell1977unstable}.  And since a schedule is given, the latenesses $\varepsilon(x)$ can themselves be idealized as a Brownian motion process.  The variance of $\varepsilon(1)$ is the variance rate of the Brownian motion, and is assumed given. It captures all random components of travel time, such as traffic and passenger boardings, and will be denoted $\sigma_0^2$. The expectation of $\varepsilon(1)$ (or ``drift''), denoted $m$, will depend on the speed of the schedule - the faster the schedule, the larger the drift.

To see how, let us introduce the average bus paces under NTP and TSP as $T_u$ and $T_c$, respectively, and denote their difference by $\Delta \equiv T_u - T_c > 0$  (Given our chosen units, pace is simply the travel time between consecutive signals and $\Delta$ is the average travel time savings afforded by TSP.) We also define the schedule pace as $T_s$.   Under NTP the drift is $m \equiv T_s - T_u$.  Thus, the expected deviation must grow without bound, unless we choose $T_s$ to match $T_u$. This is a tall order, but even if we could fulfill it the variance of the lateness would still grow without bound.  The situation is equally dire for TSP, since the same can be said about it if we just replace $T_u$ by $T_c$ in the above logic. 

Fortunately, CSP turns out to bound lateness for the broad range of schedule paces: $T_s \in (T_c, T_u)$. This comes about because with CSP the drift $m$ is state-dependent, and when the schedule pace is in $(T_c, T_u)$ the drift directions always point toward the state $\varepsilon(x) = \delta$. This can be seen from \eqref{4b}: when $\varepsilon(x) > \delta$ the bus requests priority and:

\begin{equation}
m_+ = T_c - T_s < 0;
\tag{5a}
\end{equation}

\noindent so the drift pulls large deviations down -- the plus subscript signifies that this is the drift when the lateness is greater than $\delta$. And the opposite happens when $\varepsilon(x) \le \delta$ as then the drift is:

\begin{equation}
m_- = T_u - T_s > 0. 
\tag{5b}
\end{equation}
\setcounter{equation}{5}

\begin{figure}[h]
\centering
\includegraphics[width=0.8\textwidth,trim={4cm 7cm 5.5cm 7cm},clip]{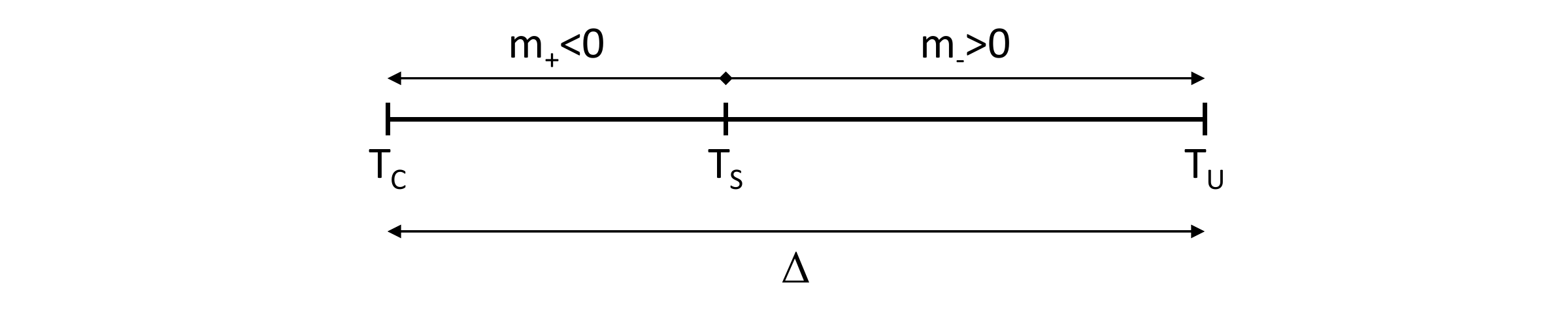}
\caption{Depiction of pace and drift}
\label{pace_figure}
\end{figure}

\noindent The relationships between pace and drift terms are depicted in \autoref{pace_figure}. Clearly, the state $\varepsilon(x) = \delta$ is a stable equilibrium. 

Buses, however, experience many perturbations that throw them off this equilibrium. To estimate the consequences, it is useful first to analyze how far a bus must travel before returning to equilibrium after an initial perturbation. To do this, call the initial deviation $\varepsilon_0 - \delta$, where $\varepsilon_0$ is the initial lateness. Now that the bus is late (or early), CSP responds by requesting priority (or denying priority) at every signal until the bus has returned to its stable equilibrium. The distance $x_0$ to achieve full recovery, i.e.: 

\begin{equation*}
x_0 \equiv \inf\{ x > 0 \vert \varepsilon(x)=\delta\},
\end{equation*}

\noindent is a first passage problem of Brownian motion theory. It is well known that $x_0$ can be represented by an inverse Gaussian distribution with parameters $\mu$ and $\lambda$; i.e.: 

\begin{equation*}
x_0 \sim IG \left( \mu, \lambda \right).
\end{equation*}

\noindent The parameters $\mu$ and $\lambda$ are functions of $(\varepsilon_0-\delta)$, $\sigma_0^2$, and the drift $m$. Because the drift is $m_+$ if $\varepsilon(x)>\delta$ and $m_-$ if $\varepsilon(x)\le\delta$, we shall abbreviate it as $m_*$. The expressions in question are then: $\mu=\frac{\varepsilon_0-\delta}{-m_*}$ and $\lambda=\frac{(\varepsilon_0-\delta)^2}{\sigma_0^2}$. Formulas for the expectation and variance of $x_0$ can then be obtained from the properties of the inverse Gaussian distribution:

\begin{equation}\label{(6)}
E[x_0]=\mu=\frac{\varepsilon_0 - \delta}{-m_*} > 0; \hskip5mm var[x_0]=\frac{\mu^3}{\lambda}=\frac{(\varepsilon_0 - \delta) \sigma_0^2}{-m_*^3} > 0.
\end{equation}

Note that the recovery rate, i.e. the deviation time recovered per unit of distance traveled as expressed by the ratio of $\vert \varepsilon_0 - \delta \vert$ to $E[x_0]$, is $\vert m_* \vert$. Also note that the moments of $x_0$ are finite. These properties imply that the system is resilient; it always recovers and returns to equilibrium from any deviation, $\varepsilon_0 - \delta$, no matter how large.

\begin{figure}[h!]
	\centering
	\includegraphics[width=0.6\textwidth]{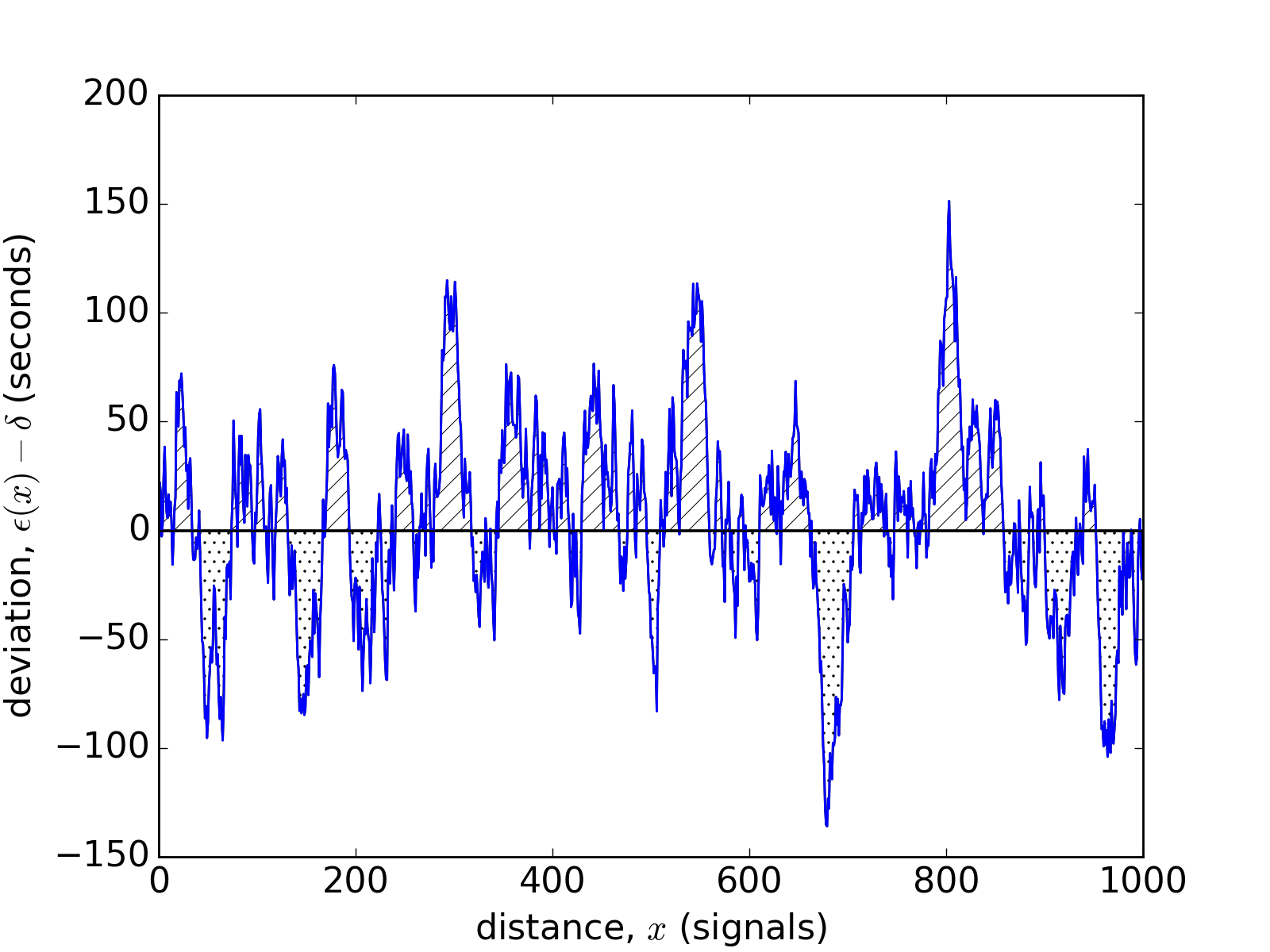}
	\caption{Deviation of a single bus with CSP}
	\label{deviation_conditional}
\end{figure}

Let us now examine how the system behaves in the long term. \autoref{deviation_conditional} shows a typical realization of deviation from equilibrium. Note how the deviations alternate in sign over consecutive time intervals and how the process spends more time on one side than the other. We shall first obtain as a preceding step the fraction of time $P_+$ that the process is on the positive side. To do this, we assume (logically) that the distribution of the initial deviation after a state $(\varepsilon_0 - \delta)$ is an even function; i.e. that positive and negative deviations are of similar character. This is reasonable given the nature of traffic and signal delays. However, since the distance to recovery $x_0$ depends on the appropriate drift value, the bus can be expected to spend more of the route on one side of equilibrium than the other unless the drifts are equal in magnitude, $\vert m_+ \vert=\vert m_- \vert$. Imagine an experiment where the bus starts at $\varepsilon_0-\delta>0$. Wait for recovery, then apply a negative deviation of the same magnitude and wait for recovery. From \eqref{(6)}, the expected distance to recovery from the positive deviation is $\frac{\varepsilon_0 - \delta}{-m_+}$ and the expected distance to recovery from the negative deviation is $\frac{-(\varepsilon_0 - \delta)}{-m_-}$. The total distance covered in this experiment is the sum of the two, $\frac{\varepsilon_0 - \delta}{-m_+}+\frac{-(\varepsilon_0 - \delta)}{-m_-}$. The fraction of distance where $\varepsilon(x) > \delta$ is therefore:

\begin{equation}\label{late}
P_+\equiv \text{fraction of distance where CSP is active} = \frac{\frac{\varepsilon_0 - \delta}{-m_+}}{\frac{\varepsilon_0 - \delta}{-m_+}+\frac{-(\varepsilon_0 - \delta)}{-m_-}}=\frac{\frac{1}{-m_+}}{\frac{1}{-m_+} + \frac{1}{m_-}}=\frac{1}{1+\gamma},
\end{equation}

\noindent where: 

\begin{equation}\label{gamma}
\gamma \equiv |m_+ / m_-| = (T_s - T_c)/(T_u - T_s) >0.
\end{equation}

\noindent Note that the parameters $\gamma=\frac{-m_+}{m_-}$ and $\Delta=m_- - m_+$ are in a 1:1 correspondence with the drifts. The inverse formulas are:

\begin{equation}\label{correspondence}
\begin{split}
m_- &= \frac{\Delta}{1+\gamma} \\
m_+ &= -\Delta \frac{\gamma}{1+\gamma} \\
\end{split}
\end{equation}

\noindent Note as well that \eqref{late} is independent of $(\varepsilon_0 - \delta)$; it only hinges on the equal magnitude of the initial deviation. Therefore, the result holds for any even distribution of initial deviation. Clearly, the fraction of distance where $\varepsilon(x) < \delta$ (CSP is dormant) is the complement of \eqref{late}:

\begin{equation}\label{early}
P_-\equiv \text{fraction of distance where CSP is dormant} =\frac{\gamma}{1+\gamma}.
\end{equation}

\noindent Notice that both fractions are $\frac{1}{2}$ if $-m_+=m_-$, as then $\gamma=1$.

We are now ready to estimate the mean and variance of the steady state deviation. We will do this by estimating separately the moments of the positive and negative parts of the deviation, and then using \eqref{late} and \eqref{early} to combine. Now picture the trajectory of a bus that passes through the equilibrium state many times, such as the one shown in \autoref{deviation_conditional}. We want to estimate $E[(\varepsilon-\delta)^+]$ and $E[((\varepsilon-\delta)^+)^2]$ considering only the shaded areas. And do the same for the dotted areas to estimate $E[(\varepsilon-\delta)^-]$ and $E[((\varepsilon-\delta)^-)^2]$.

In the shaded areas the process behaves as if it had a reflecting barrier at $(\varepsilon-\delta)^+=0$ and a negative drift. The equilibrium solution for this case is that $(\varepsilon-\delta)^+$ is exponentially distributed with scale parameter $\beta_+=\frac{\sigma_0^2}{-2m_+}$ \citep{newell1977unstable}. The first and second moments can then be obtained from the properties of the exponential distribution:

\begin{equation}\label{epsilon+}
E[(\varepsilon-\delta)^+]=\beta_+=\frac{\sigma_0^2}{-2m_+}; \hskip5mm E[((\varepsilon-\delta)^+)^2]=2\beta_+^2=\frac{\sigma_0^4}{2m_+^2}.
\end{equation}

\noindent Similarly, in the dotted areas the process has a reflecting barrier at $(\varepsilon-\delta)^-=0$ and positive drift. The negative deviation is also exponentially distributed, but with scale parameter $\beta_-=\frac{\sigma_0^2}{2m_-}$. The first and second moments are:

\begin{equation}\label{epsilon-}
E[(\varepsilon-\delta)^-]=\frac{\sigma_0^2}{2m_-}; \hskip5mm E[((\varepsilon-\delta)^-)^2]=\frac{\sigma_0^4}{2m_-^2}.
\end{equation}

We shall then combine these estimates with \eqref{late} and \eqref{early}, $P_+$ and $P_-$, to obtain:

\begin{equation}\label{4.2}
\begin{split}
E[\varepsilon-\delta]&=E[(\varepsilon-\delta)^+] P_+ + E[(\varepsilon-\delta)^-] P_- \\
&=\frac{\sigma_0^2}{-2m_+}\left( \frac{1}{1+\gamma} \right) + \frac{-\sigma_0^2}{2m_-}\left( \frac{\gamma}{1+\gamma} \right)\\
&=\frac{\sigma_0^2}{2\Delta}\left( \frac{1-\gamma^2}{\gamma} \right), \\
\end{split}
\tag{13a}
\end{equation}

\noindent and:

\begin{equation}\label{varlate}
\begin{split}
var[\varepsilon-\delta]&=E[(\varepsilon-\delta)^2] - E[\varepsilon-\delta]^2 \\
&=\frac{\sigma_0^4}{2m_+^2}P_+ + \frac{\sigma_0^4}{2m_-^2}P_- - \frac{\sigma_0^4}{4\Delta^2}\left(\frac{(1+\gamma)^2(1-\gamma)^2}{\gamma^2}\right) \\
&=\frac{\sigma_0^4}{4\Delta^2}\left( \frac{(1+\gamma)^2(1+\gamma^2)}{\gamma^2} \right). \\
\end{split}
\tag{13b}
\end{equation}
\setcounter{equation}{13}

\noindent The last equality in the above is the result of eliminating $(P_+, P_-, m_+, m_-)$ using \eqref{late}, \eqref{correspondence}, and \eqref{early} and then simplifying the resulting expression.

Note $\delta$ does not appear in either expression. Therefore, for any given $T_s$, its optimum is the value that achieves $E[\varepsilon]= 0$; i.e., $\delta^*$ is the negative of the expression for $E[\varepsilon-\delta]$:

\begin{equation}\label{delta*}
\delta^*=-\frac{\sigma_0^2}{2\Delta}\left( \frac{1-\gamma^2}{\gamma} \right).
\end{equation}

\subsection{Simulation Tests}

The formulas derived in \eqref{(6)}, \eqref{4.2}, and \eqref{varlate} contain many simplifications. To check their accuracy, they were compared against a more realistic microsimulation that tracked many hundreds of buses over a very long street with both stations and pre-timed signals with arbitrary offsets and a common cycle. The simulation considers that there may also be bus routes on cross streets that use signal priority, leading to conflicting priority requests. These conflicts were resolved on a first-come first-served basis. The parameters of the street and signals are given in \autoref{sim_parameters}. 

\autoref{components} presents the components of bus travel time. The first column, $\bm{T_l}$, follows from the data in \autoref{sim_parameters} and an assumed bus cruising speed of 30 mi/hr. Line haul represents the freeflow travel time between two signals and therefore has 0 variance. The second column, $\bm{T_p}$, is a normal random variable whose parameters were chosen to represent a situation with considerable demand and variability. The third and fourth columns were calculated with standard traffic engineering methods using the data of \autoref{sim_parameters}. The derivations are too long to be included here, but the results were confirmed by the simulations. The last two columns give the moments of the combined time. The variance entries give the values of the parameter $\sigma_0^2$ appearing in (13) and \eqref{delta*}. They are calculated assuming that the components of $\bm{T_u}=\bm{T_l}+\bm{T_p}+\bm{T_n}$ and $\bm{T_c}=\bm{T_l}+\bm{T_p}+0.08\bm{T_n}+0.92\bm{T_y}$ are independent\footnote{The signal delay component of $\bm{T_c}$ is taken to be $0.08\bm{T_n}+0.92\bm{T_y}$ rather than $\bm{T_y}$ because with the configuration of \autoref{sim_parameters} about 8\% of the requests for priority are denied due to conflicting crossing buses.}.

Recall that the two control parameters are $T_s$ and $\delta$. The first simulations focus on the relationship between $T_s$ and the formulas in \eqref{(6)}, \eqref{4.2}, and \eqref{varlate}, so $\delta$ is held constant at 0. The second simulation holds $T_s$ constant at 49.21 ($\gamma=1$) and tests the formula for $\delta^*$ in \eqref{delta*}. The final simulation uses $\delta=0$ and $T_s=49.21$ to compare CSP to NTP and TSP.

\begin{table}[h]
\centering
\caption{System configuration values}
\label{sim_parameters}
\begin{tabular}{|l|r|l|}
\hline
name&value&units \\
\hline
Signal separation&0.25&mi \\
\hline
Station separation&0.25&mi \\
\hline
Crossing bus route separation&0.25&mi \\
\hline
Crossing route headway&10&min \\
\hline
Cycle length&100&s \\
\hline
Green phase [used by bus]&60&s \\
\hline
Advance notice for priority&10&s \\
\hline
Clear lag [to safely end phase]&20&s \\
\hline
\end{tabular}
\end{table}

\begin{table}[h]
\centering
\caption{Components of bus travel time}
\label{components}
\begin{tabular}{|r|>{\raggedleft}p{2cm}|>{\raggedleft}p{2cm}||>{\raggedleft}p{2cm}|>{\raggedleft\arraybackslash}p{2cm}|>{\raggedleft\arraybackslash}p{2cm}|>{\raggedleft\arraybackslash}p{2cm}|}
\hline
&	$\bm{T_l}$: line haul time&	$\bm{T_p}$: extra delay due to traffic / passengers&	$\bm{T_n}$: signal delay (no priority)&	$\bm{T_y}$: signal delay (yes priority)&	$\bm{T_u}=\bm{T_l}+\bm{T_p}+\bm{T_n}$&	$\bm{T_c}=\bm{T_l}+\bm{T_p}+0.08\bm{T_n}+0.92\bm{T_y}$\\
\hline
expectation&	$T_l=$30.00 $s$ &	$T_p=$13.60 $s$&	$T_n=$8.20 $s$&	$T_y=$2.55 $s$&	$T_u=$51.80 $s$&	$T_c=$46.62 $s$ \\
\hline
variance&	$\sigma_l^2=$0.00 $s^2$&	$\sigma_p^2=$130.9 $s^2$&	$\sigma_n^2=$154.2 $s^2$&	$\sigma_y^2=$17.35 $s^2$&	$\sigma_{0u}^2=$285.0 $s^2$&	$\sigma_{0c}^2=$159.6 $s^2$ \\
\hline
\end{tabular}
\end{table}

\begin{figure}
	\centering
	\includegraphics[width=0.6\textwidth]{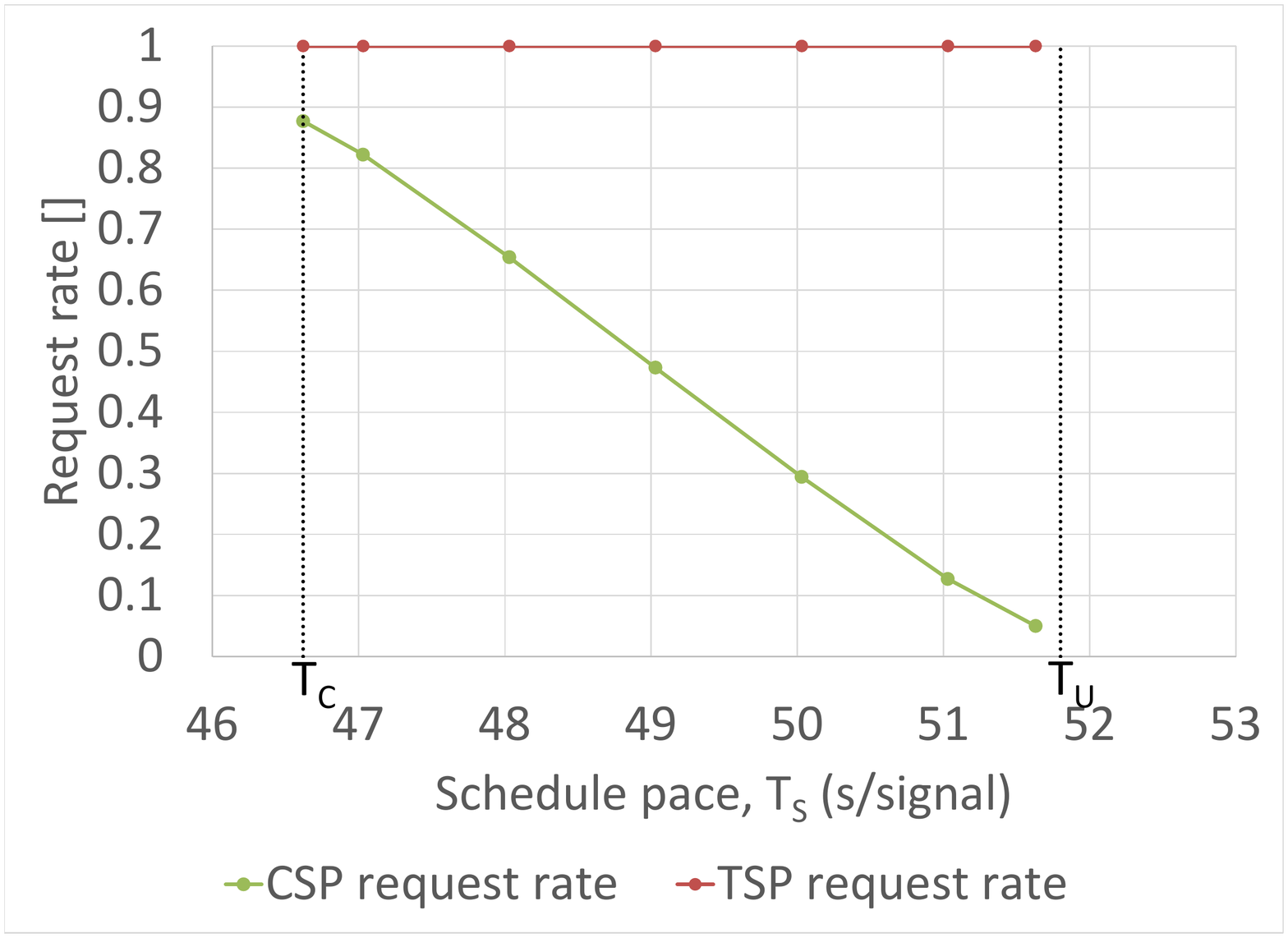}
	\caption{Signal priority request rate vs schedule pace, no holding at stations}
	\label{nohold_request}
\end{figure}

For CSP, the expectation of $\bm{T_s}$ is a control variable but its variance, $\sigma_{0s}^2$, depends on the frequency $r$ of priority requests, which is unknown. Since in practical applications this parameter needs to be known beforehand we did not evaluate it from simulation. Instead we used a calculation method that could be used in practice (and verified it with simulations). Reasonably, we assumed that $\sigma_{0s}^2$ could be linearly interpolated between the known values of $\sigma_{0c}^2$ and $\sigma_{0u}^2$ as per:

\begin{equation}\label{*}
\sigma_{0s}^2=(1-r)\sigma_{0u}^2 + r \sigma_{0c}^2,
\end{equation}

\noindent and also assumed that the frequency of priority requests equals the fraction of distance that the bus is lagging behind; i.e. that:

\begin{equation}\label{**}
r \approx P_+ = \frac{1}{1+\gamma}.
\end{equation}

Thus,

\begin{equation}\label{***}
\sigma_{0s}^2=\frac{\left( \gamma\sigma_{0u}^2 + \sigma_{0c}^2\right)}{1+\gamma}.
\end{equation}

\noindent \autoref{nohold_request} confirms the reasonableness of \eqref{**}.

\begin{figure}
	\centering
	\begin{subfigure}[b]{0.43\textwidth}
		\centering
		\includegraphics[width=\textwidth]{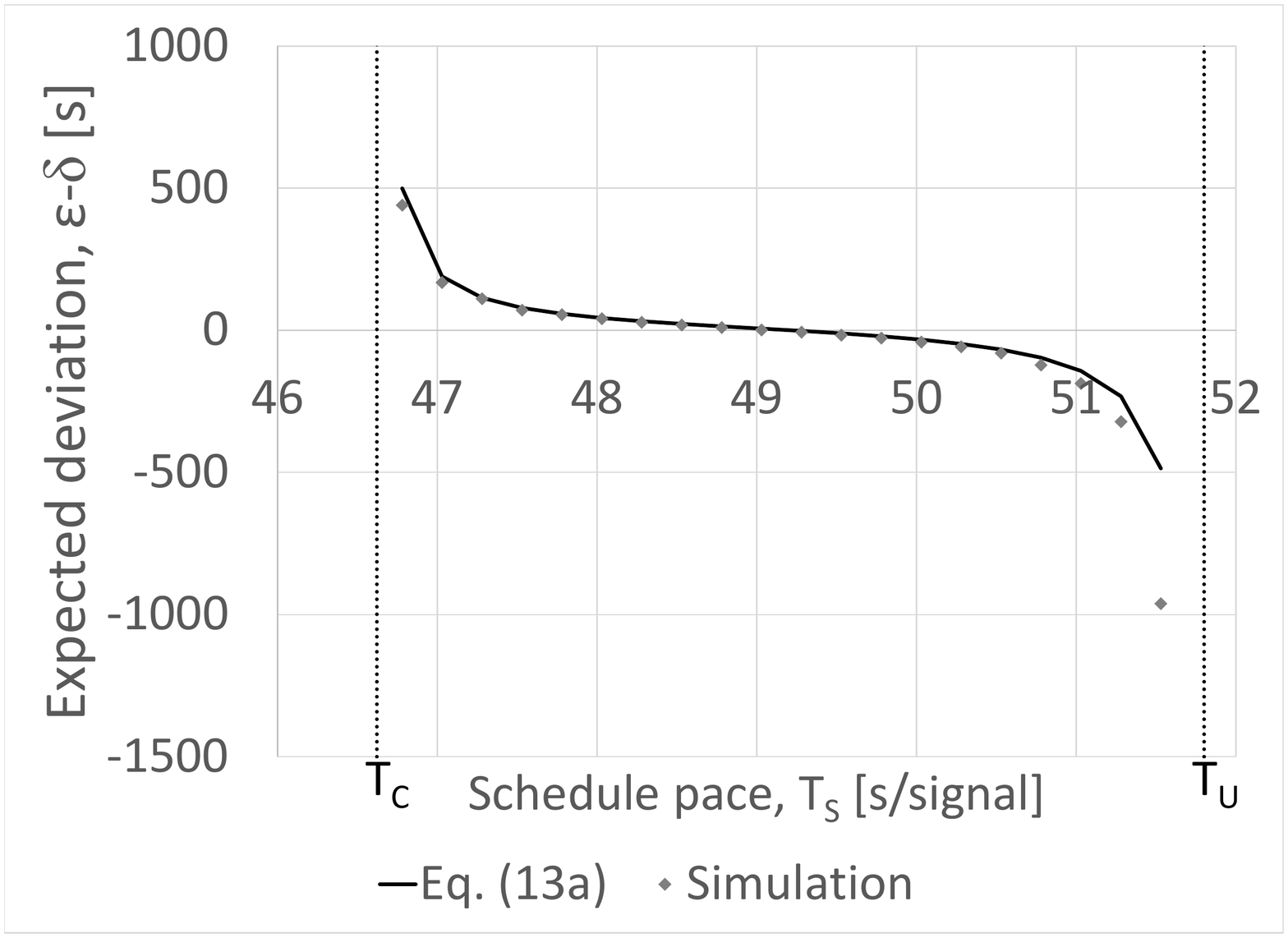}
		\caption{Predicted and actual expected lateness for different schedule paces, $\delta=0$}
		\label{epsilon}
	\end{subfigure}
	~
	\begin{subfigure}[b]{0.43\textwidth}
		\centering
		\includegraphics[width=\textwidth]{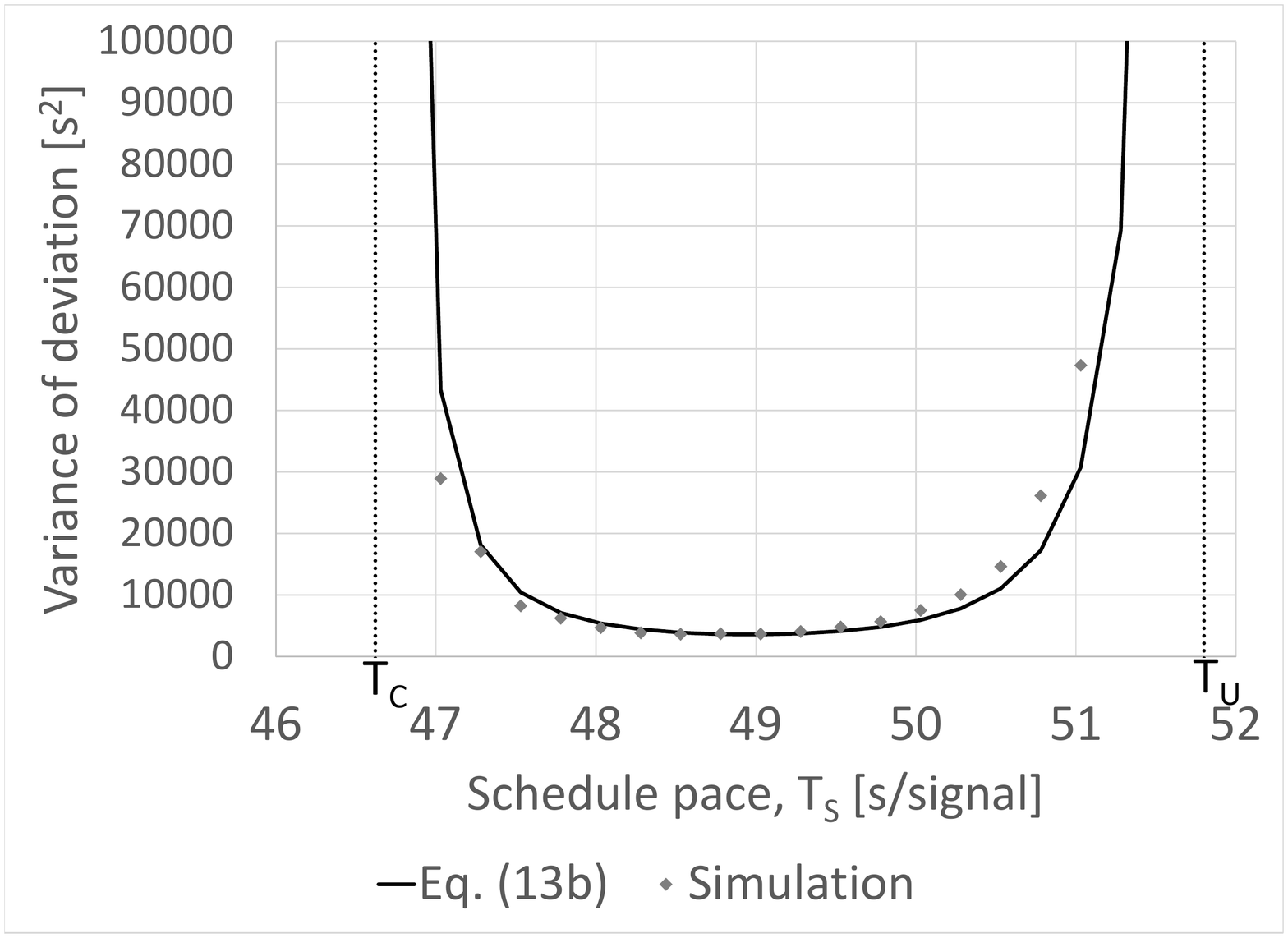}
		\caption{Predicted and actual variance of lateness for different schedule paces, $\delta=0$}
		\label{Brownian_variance}
	\end{subfigure}
	
	\begin{subfigure}[b]{0.43\textwidth}
		\centering
		\includegraphics[width=\textwidth]{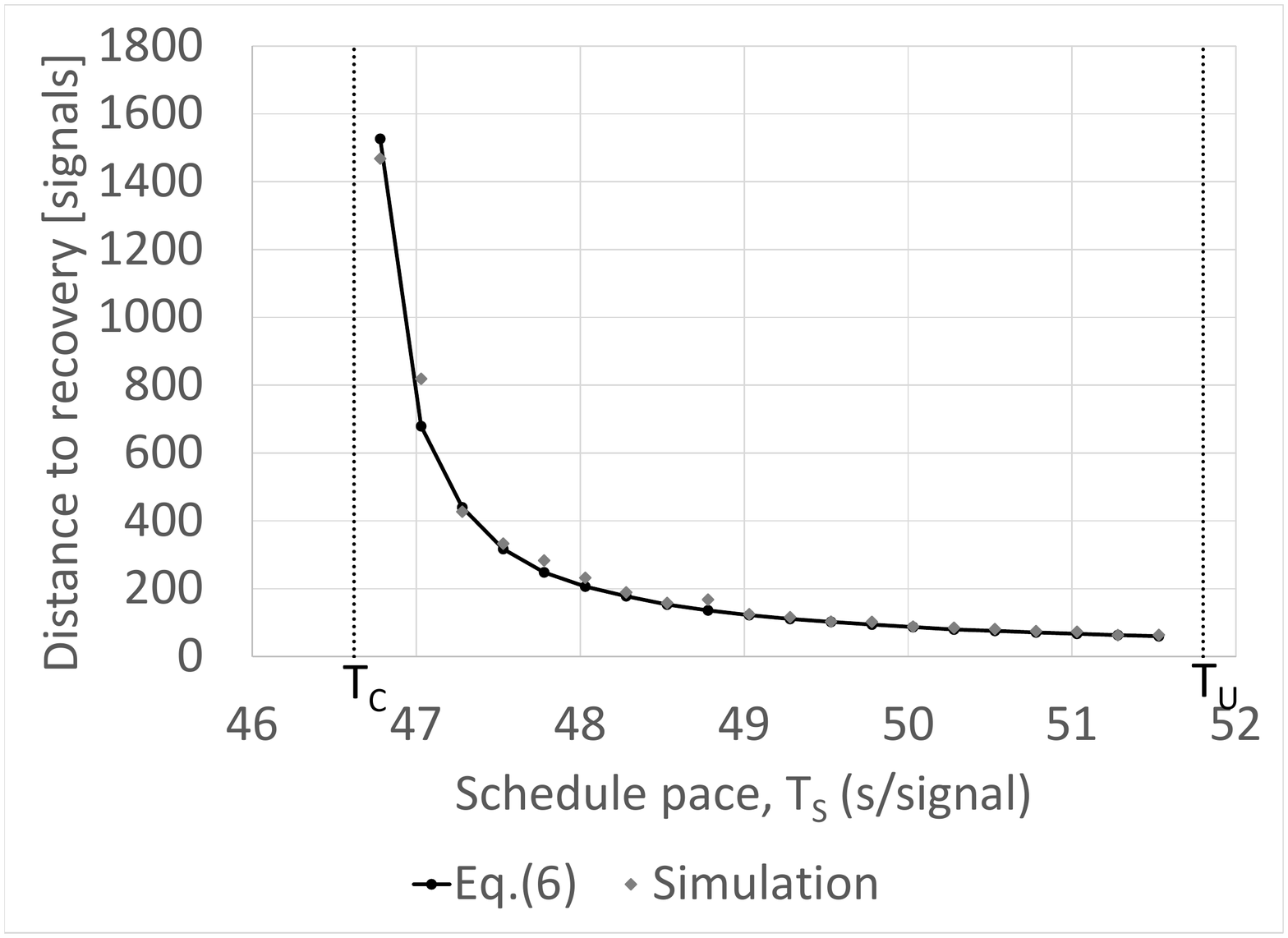}
		\caption{Predicted and actual expected times to recovery for different schedule paces, $\delta=0$}
		\label{Brownian_delay_mean}
	\end{subfigure}
\caption{Accuracy of Brownian motion formulas for three CSP metrics.}
\label{verification}
\end{figure}

Using the values of \autoref{components} and \eqref{***} the expressions for the expectation and variance of lateness in \eqref{4.2} and expected distance to recovery in \eqref{(6)} were plotted vs. $T_s$. The results are the solid lines in \autoref{epsilon}, \autoref{Brownian_variance}, and \autoref{Brownian_delay_mean}, respectively.

The accuracy of \eqref{4.2} and \eqref{(6)} was then assessed by comparing simulated outcomes to the predictions. The simulated outcomes are shown as dots on \autoref{epsilon}, \autoref{Brownian_variance}, and \autoref{Brownian_delay_mean}. Note the good match on all three plots.

\subsection{Interpretation and Fine-Tuning}

\begin{figure}
	\centering
	\begin{subfigure}[b]{0.43\textwidth}
		\centering
		\includegraphics[width=\textwidth]{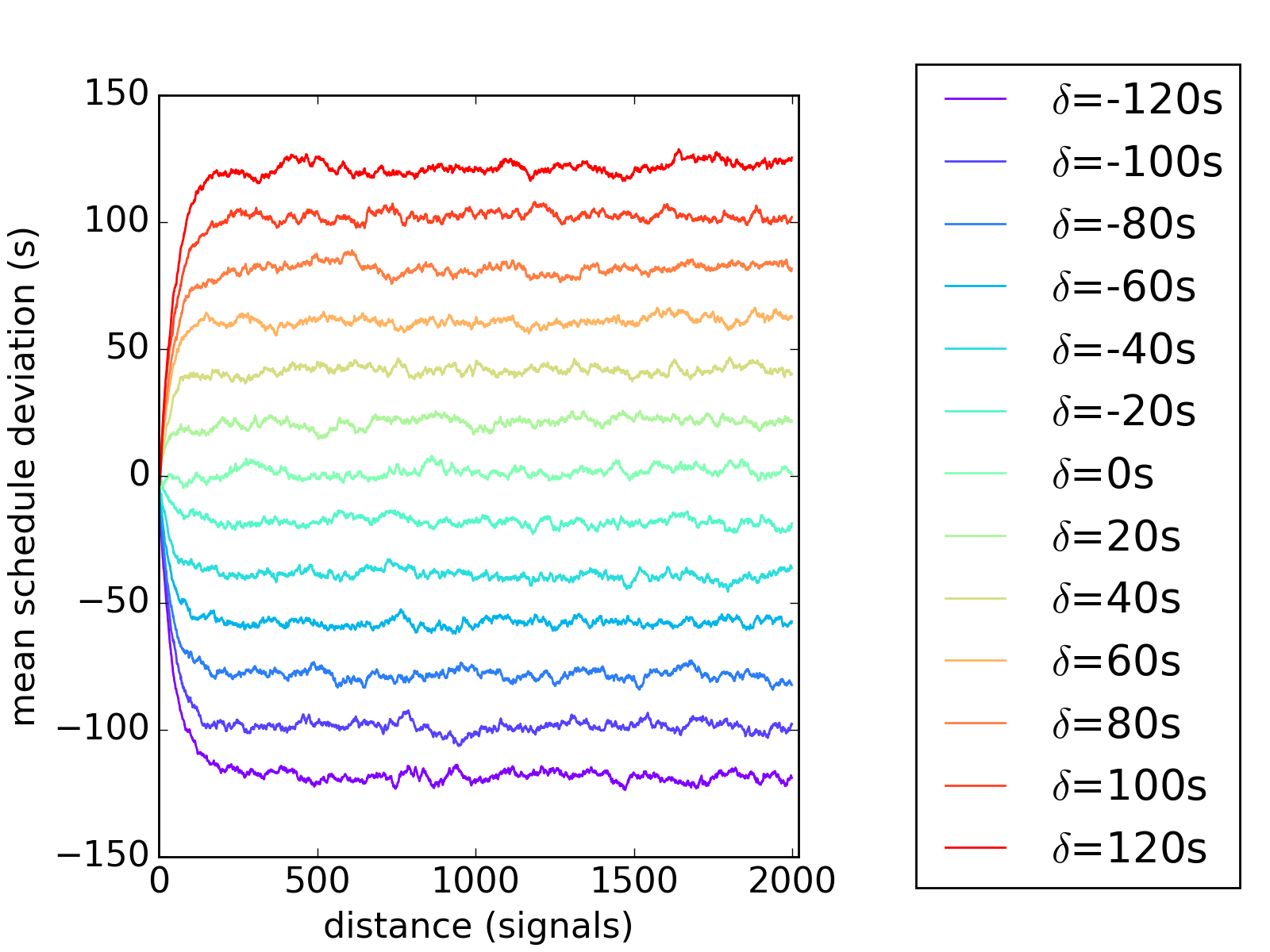}
		\caption{Sample paths of lateness under CSP with different $\delta$'s: case with $\gamma = 1$.}
		\label{thresholds}
	\end{subfigure}
	~
	\begin{subfigure}[b]{0.43\textwidth}
		\centering
		\includegraphics[width=\textwidth]{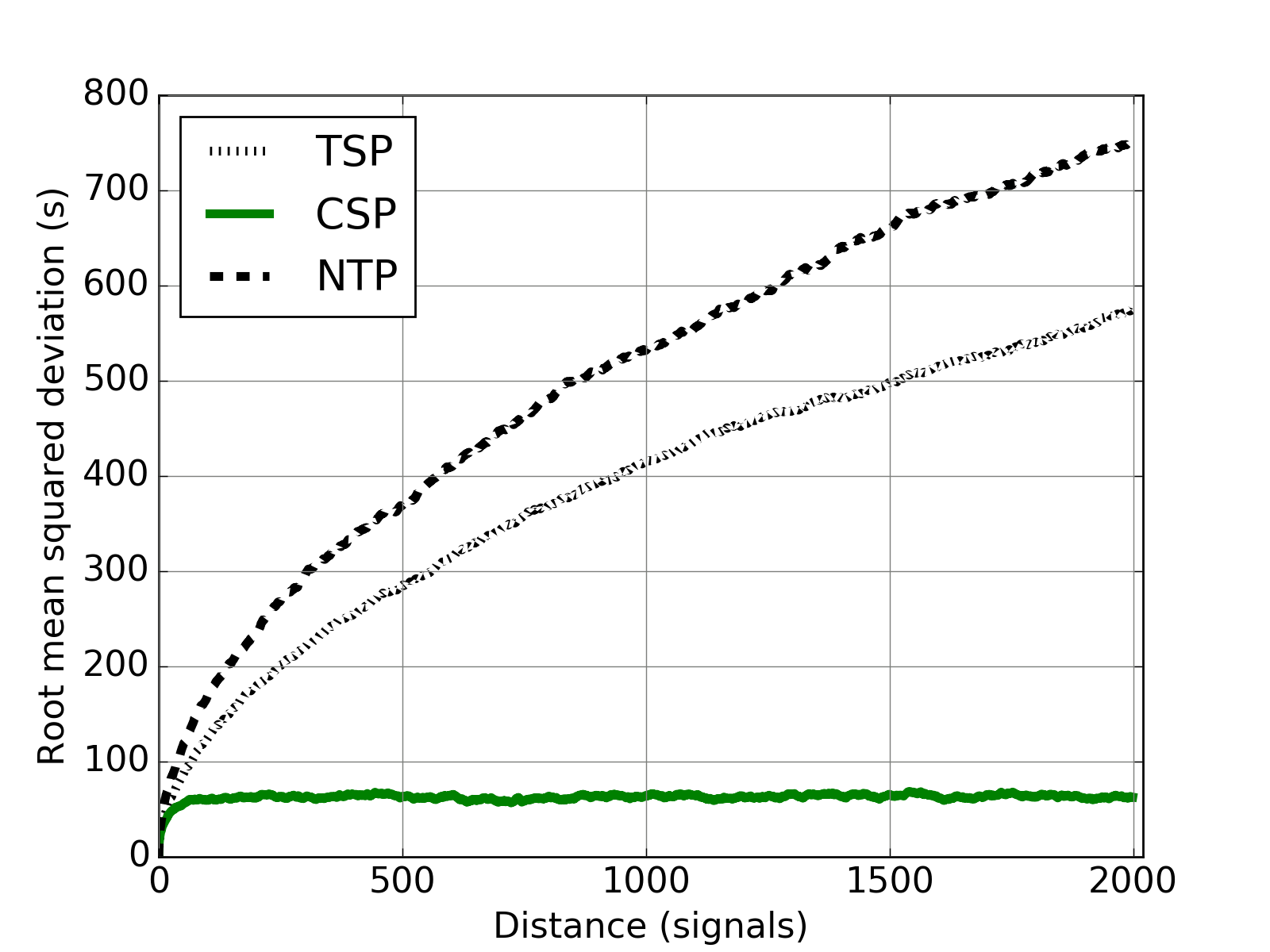}
		\caption{RMS of schedule deviation vs distance for three strategies: optimized parameters.}
		\label{paxperspective}
	\end{subfigure}
\caption{Sample results for case of long headways with schedule}
\label{sample}
\end{figure}

The simulations also confirm that CSP can be very effective if $\delta$ is chosen properly, and that it vastly outperforms NTP and TSP. The first idea is illustrated by \autoref{thresholds}, which considers a schedule with $\gamma = 1$ ($T_s = 49.21s$) and shows sample paths of lateness for different $\delta$'s. Note the importance of choosing $\delta$ properly, and how the optimum value arising from \eqref{4.2}, $\delta^* = 0$, performs best.  The second idea is illustrated by \autoref{paxperspective}. It displays the root mean squared (RMS) schedule deviations computed over 10,000 draws for CSP, NTP and TSP, using for each of these three strategies the most favorable parameters; i.e., $\gamma =0$ for NTP, $\gamma = \infty$ for TSP and $\gamma = 1, \delta = 0$ for CSP.  

Choosing a schedule pace, $T_s$, involves a tradeoff between the speed of the schedule and the variance of the lateness. \autoref{Brownian_variance} shows that the variance of the lateness improves with increasing $T_s$ in the interval $[T_c, \frac{T_u+T_c}{2}]$. Clearly thus, $T_s^*$ must be in this interval. For example, if the bus operator wishes to constrain the variance of the steady state deviation to some upper limit $\sigma_m^2$, the optimum $T_s$ is:

\begin{equation}\label{tstar}
T_s^*=\min\{T_s: var[\varepsilon-\delta] \leq \sigma_m^2\},
\end{equation}

\noindent which must occur where the constraint is binding. Therefore, in this case $T_s^*$ occurs where the declining branch of the curve in \autoref{Brownian_variance} reaches the value $\sigma_m^2$.

\section{Case (ii): Holding by schedule}

In Case (ii), both drivers and signals act as control agents. We assume once again that there is a schedule, but now passengers consult it before arriving at stations to choose their arrival times. Drivers check the schedule at every station (i.e. every station is a control point) and hold until the scheduled departure time if they are early. With this form of operation the number of passengers picked up by a bus is independent of the lateness. Thus, there is no bunching tendency. As a result buses will be treated in isolation as in Section 4.

The following subsection, 5.1, explains the modifications to the Brownian Motion model of Subsection 4.1 that are needed. Subsection 5.2 tests the main result with simulations. Subsection 5.3 shows how to optimize the control parameters and discusses the results.

\subsection{Brownian Motion}
The Brownian motion model developed in Case (i) is now extended to Case (ii). The act of holding by schedule is incorporated into the model by turning the equilibrium state, represented by the line $\varepsilon(x)=\delta$, into a reflecting barrier. This is an approximation because lateness is continuous and holding is applied only at stations, but passengers can only board and alight at stations and the approximation error is small when every station is a control point. In Case (i) with CSP, formulas for the expectation and variance of steady state deviation, \eqref{4.2} and \eqref{varlate}, were obtained by deriving expressions for the positive and negative parts of the steady state deviation, \eqref{epsilon+} and \eqref{epsilon-}, and then weighting by the fraction of distance that the bus spends in each state, \eqref{late} and \eqref{early}. In Case (ii), the reflecting barrier means that buses are never early ($m_-=\infty$, $P_-=0$). Therefore, the mean steady state deviation with CSP is simply the positive part of steady state deviation from \eqref{epsilon+}:

\begin{equation}\label{caseiiexpectation}
E[\varepsilon-\delta]=\frac{\sigma_0^2}{-2m_+}, \hskip5mm (m_+ <0)
\end{equation}

\noindent and the variance can be calculated using the first and second moments of the positive deviation from \eqref{epsilon+}:

\begin{equation}\label{caseiivar}
\begin{split}
var[\varepsilon-\delta]&=\frac{\sigma_0^4}{2m_+^2} - \left( \frac{\sigma_0^2}{-2m_+} \right)^2 \\
&=\frac{\sigma_0^4}{2m_+^2} - \frac{\sigma_0^4}{4m_+^2} \\
&=\frac{\sigma_0^4}{4m_+^2}. \\
\end{split}
\end{equation}

Because passengers know the schedule, their waiting time is proportional to the root mean square of the latenesses, $RMS[\varepsilon]$, which based on \eqref{caseiiexpectation} and \eqref{caseiivar} is:

\begin{equation}
RMS[\varepsilon]=\sqrt{E[\varepsilon]^2+var[\varepsilon]}=\sqrt{\left(\frac{\sigma_0^2}{-2m_+}+\delta\right)^2 + \frac{\sigma_0^4}{4m_+^2}}.
\end{equation}

This will be our measure of performance, assuming that $\delta$ has been optimized to minimize $RMS[\varepsilon]$; i.e. that:

\begin{equation}
\delta^*=\frac{\sigma_0^2}{2m_+}.
\end{equation}

Thus our measure of performance is:

\begin{equation}
RMS[\varepsilon]=\sqrt{var[\varepsilon-\delta]}=-\frac{\sigma_0^2}{2m_+}.
\end{equation}

The above applies to CSP. Buses with TSP also use holding by schedule. Therefore, the above expressions also apply to buses with TSP. For buses with NTP, the system is stable only if $T_s > T_u$. The above expressions apply to buses with NTP provided that $m_+$ is replaced with $-m\equiv T_u-T_s<0$. Otherwise $E[\varepsilon-\delta]=\infty$. 

\subsection{Simulation Tests}
\begin{figure}[h!]
	\centering
	\includegraphics[width=0.6\textwidth]{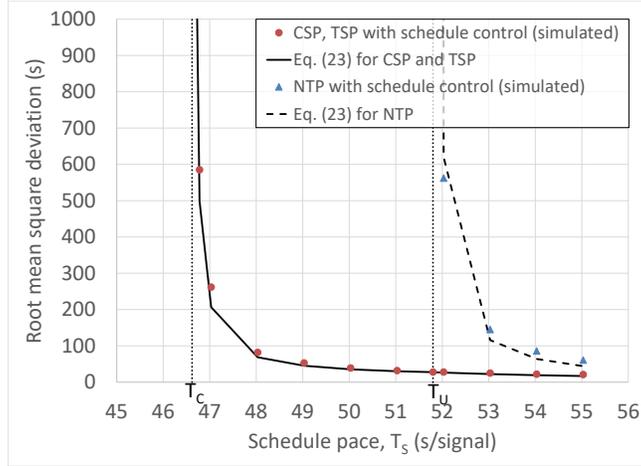}
	\caption{RMS of schedule deviation vs schedule pace}
	\label{schedule_RMSdev}
\end{figure}

The validity of the formulas for Case (ii) is evaluated using the same simulation from Case (i). \autoref{schedule_RMSdev} shows the predicted and simulated RMS schedule deviations for the CSP and TSP scenarios and a range of schedule paces. The labels $T_c$ and $T_u$ on the x-axis denote the average pace with priority at all signals and no transit priority, respectively. 

\subsection{Interpretation and Fine-Tuning}

\begin{figure}
	\centering
	\includegraphics[width=0.6\textwidth]{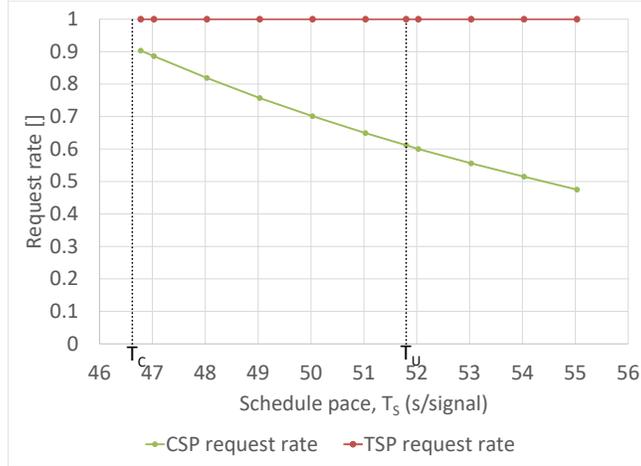}
	\caption{Signal priority request rate vs schedule pace, holding by schedule}
	\label{schedule_request}
\end{figure}

First, note that TSP and CSP have the same effect on bus performance since in both scenarios buses that are behind schedule request priority and buses that are ahead of schedule are held at the next station. The difference between TSP and CSP, then, is in their impact on other traffic. \autoref{schedule_request} shows the priority request rate and for both strategies and a range of schedule paces. The relationship between schedule pace and the priority request rate is observed to be approximately linear. It can be seen that CSP significantly reduces the request rate compared with TSP. Note there is a significant improvement even if the schedule is the fastest possible.

Comparing \autoref{schedule_request} with \autoref{nohold_request} we see that the request rates with CSP are higher in Case (ii) than Case (i). At $T_s=48$ $s$ per signal, the request rate is approximately 0.82 in Case (ii) and 0.65 in Case (i). This should be expected because by holding buses at stations, Case (ii) buses are less likely to be ahead of schedule and therefore be more likely to request priority.

As in Case (i), the choice of schedule pace, $T_s$ is a tradeoff between schedule speed and the variance of steady state deviation (or RMS deviation). As can be seen from \autoref{schedule_RMSdev}, holding by schedule eliminates the increase in variance seen in \autoref{Brownian_variance} as $T_s \rightarrow T_u$. This change means that $T_s^*$ can now be larger than $\frac{T_u+T_c}{2}$ and yield a lower minimum\footnote{The reader can verify that $\lim_{T_s \to \infty} RMS[\varepsilon]=0$.}. If desired we can still use \eqref{tstar} to determine $T_s$.

\section{Case (iii): Holding by headways without a schedule}

The third case disposes with the schedule. Without it passengers arrive at the station when they are ready to depart. This passenger behavior produces a feedback loop which makes high-frequency bus operation unstable: buses following long headways tend to encounter more passengers than usual, and the opposite occurs to buses following short headways. Because the dwell time at stations depends on the number of passengers encountered, over time the long headways get longer and the short headways tend to 0. This phenomenon is known as ``bus bunching''.

Bus bunching can be alleviated by non-linear control strategies that consider multiple buses \citep{daganzo2011reducing, xuan2011dynamic, bartholdi2012self, daganzo2017unpublished}. The ``robust control'' from \citep{daganzo2017unpublished} will be used in our evaluation because it is resistant to large disturbances. Using the notation of \eqref{one} and \eqref{two}, the holding time for this strategy is:

\begin{equation}\label{holding}
\begin{split}
D_{n,s}&=k_0 - k_1 (a_{n,x}-a_{n-1,x}) + k_1(\alpha_{n+1,x}-a_{n,x}), \\
H_{n,s}&=[D_{n,s}]^+, \\
\end{split}
\end{equation}

\noindent where $k_0$ and $k_1$ are control parameters. The $k_0$ parameter is supposed to be non-negative and is set to 0 in order to minimize holding times. The $k_1$ parameter governs the response to perturbations and a value of 0.2 works well for a broad range of demand conditions including those that will be used. This driver holding is the only strategy used in the NTP scenario. In the CSP scenario, priority is requested as per \eqref{4a} and we will use $\delta=0$. In the TSP scenario, priority is requested at all signals.

\subsection{Results and Interpretation}

The mathematical model developed in Cases (i) and (ii) is not applicable to Case (iii) because buses can no longer be treated in isolation. The interactions between buses are complex and nonlinear as \eqref{holding} shows. Every bus is both ahead of and behind every other bus and perturbations to any one bus can propagate to all others. Because of these complexities, we will consider Case (iii) only in simulation. 

Since no schedule is used it will be assumed that all passengers arrive randomly, that their arrival times follow a continuous uniform distribution, and that each boarding passenger delays its bus by a set amount of time. To avoid boundary effects from the first and last bus, we choose to model Case (iii) as a closed loop. This loop could represent a bus route where service is provided by a dedicated fleet of buses. The simulated loop consists of 40 segments, each containing one signal and one station. The passenger demand rate is 0.9375 per minute per station and the boarding delay per passenger 2.0 $s$. All segments are the same and the dimensions are identical to those used for Cases (i) and (ii). The simulation is initialized with a specified number of buses placed at random locations on the loop and runs for a duration of 10 hours. Each signal control scenario in Case (iii) is simulated for a range of fleet sizes.

\begin{figure}
	\centering
	\begin{subfigure}[b]{0.8\textwidth}
		\centering
		\includegraphics[width=\textwidth]{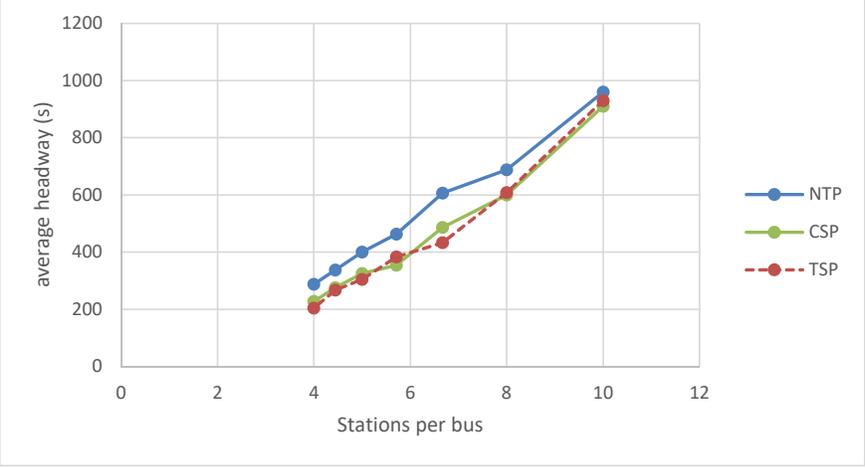}
		\caption{Average headway vs stations per bus}
		\label{robust_avg_headway}
	\end{subfigure}
	~
	\begin{subfigure}[b]{0.8\textwidth}
		\centering
		\includegraphics[width=\textwidth]{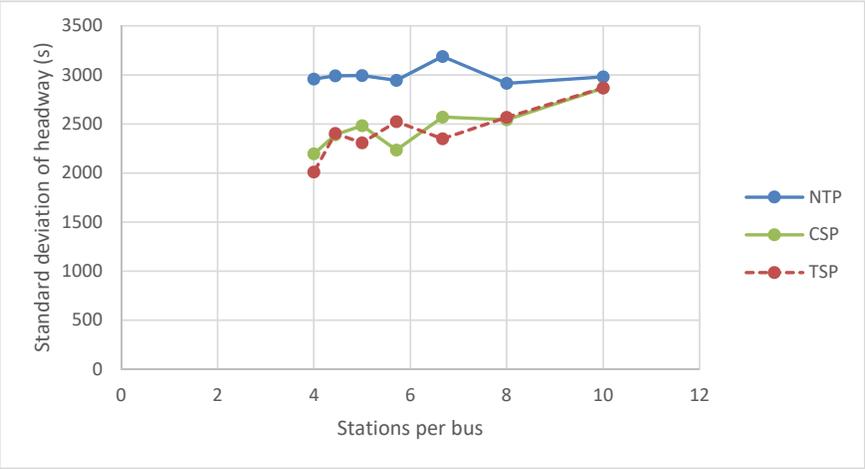}
		\caption{Standard deviation of headway vs stations per bus}
		\label{robust_stdev_headway}
	\end{subfigure}
\caption{Average and standard deviation of headway, holding by headways}
\label{avgstdheadway}
\end{figure}

\autoref{avgstdheadway} shows the average and standard deviation of headway for the three scenarios and for different numbers of buses, as represented by the number of stations per bus. As in Case (ii), the use of driver holding means that CSP and TSP have the same effect on bus operations. In both cases, the effect of signal priority is to reduce both the average headway and the standard deviation of headway. This is good. Reductions in the average headway imply that buses travel faster, since each simulation uses a fixed number of buses. Therefore, lower average headways have multifaceted benefits: they imply less waiting, shorter rides, and higher bus productivity. The decrease in the standard deviation reflects an improvement in reliability. 

\begin{figure}
	\centering
	\includegraphics[width=0.8\textwidth]{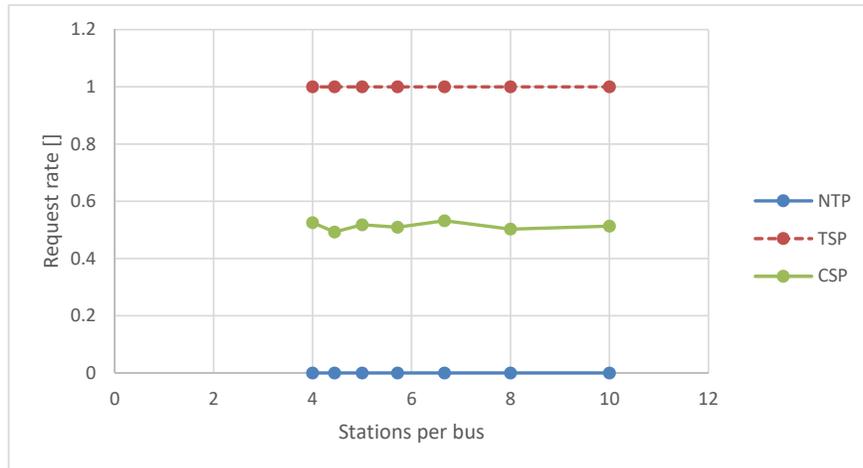}
	\caption{Signal priority request rate vs schedule pace, holding by headways}
	\label{robust_request}
\end{figure}

Although CSP and TSP have the same effect on bus operations, they differ in their impact on other traffic. \autoref{robust_request} shows the rate at which priority requests are sent for the three scenarios and the same range of bus densities as \autoref{avgstdheadway}. The request rates of NTP and TSP are fixed, and the rate with CSP is stable at about 0.5 regardless of bus density. Since CSP sends only half as many priority requests as TSP, it will have less impact on other traffic. This is its main advantage as compared to TSP.

\section{Conclusion}

In current practice, buses often struggle to stay on schedule and/or to maintain regular headways. Transit agencies design schedules with considerable slack so that most buses arrive early and instruct drivers to hold at stations until the scheduled departure time. These practices can smooth out small perturbations but do not fully address the reliability problems of buses and significantly increase their average pace. At the same time, transit signal priority is used in many places to reduce the signal delay that buses experience, but the potential for signals to serve as an additional control agent working on bus reliability has not been systematically evaluated. This paper fills the gap by considering 9 scenarios representing the combinations of three forms of driver control \{none, holding by schedule, holding by headways\} and three forms of signal priority \{no transit priority (NTP), conditional signal priority (CSP), transit signal priority (TSP)\}. The analysis is organized into three cases, each one covering one form of driver control.

Case (i), with no holding at stations, is used by some infrequent bus services. Since headways are long and passenger demand is small, there is little interaction between buses and each one can be analyzed independently. A mathematical model based on Brownian motion is proposed and is used to develop formulas for the expectation and variance of the steady state deviations. Both NTP and TSP can have an expectation of 0 with the appropriate schedule, but their variances tend to infinity. CSP has an expectation of 0 and a minimum (finite) variance with a schedule pace halfway between NTP and TSP. These findings are verified in simulation. Schedule pace is also found to have a linear relationship with the request rate of signal priority in the CSP scenario.

Case (ii), with holding by schedule, is the most common way of operating bus service today. The Brownian motion model from Case (i) is also applicable to this case. Holding is applied at every station, which turns the equilibrium state into a reflecting barrier (i.e. earliness is corrected immediately). NTP, CSP, and TSP can all achieve finite variance of steady state deviation given sufficient slack. CSP and TSP now perform equally well for bus operations, but CSP has less impact on other traffic because it sends fewer priority requests. Because earliness is corrected by holding, the request rate of signal priority with CSP is higher than in Case (i). The relationship between schedule pace and the request rate remains linear.

Case (iii), with holding by headways, has been proposed for high-frequency bus service but is not in widespread use. Because passengers arrive randomly with respect to the schedule, buses are affected by the position of other buses and the Brownian motion model of Case (i) is not applicable. A closed-loop simulation is used instead. As in Case (ii), CSP and TSP are found to have the same impact on bus operations. Signal priority reduces the average headway by making buses go faster and decreases the standard deviation of headway by allowing the system to recover more quickly from perturbations. The request rate of signal priority with CSP is around 0.5 regardless of bus density in the loop. The disruption caused by priority requests is currently a barrier to implementing TSP on high-frequency bus routes. By cutting the number of requests in half, CSP would allow for implementation in more cases.

The analytical and simulation results presented here demonstrate that CSP can be used not just to accelerate buses but to improve their reliability in all types of scenarios. It can do so as well as TSP in some scenarios and better in others. And in all scenarios is has less impact on traffic by sending fewer priority requests. Future work will consider practical constraints such as integration of signal priority with signal progression for cars and variable block geometry along a corridor.

\bibliographystyle{plain}
\bibliography{my_bibliography}

\end{document}